\documentstyle{amsppt}
\input amssym
\loadbold
\def\lh{\left\{}
\def\rh{\right\}}
\topmatter
\title
Infinite $J$-matrices and a matrix moment problem
\endtitle
\author
M. Krein
\endauthor
\thanks
Dokl.\ Akad.\ Nauk SSSR {\bf 69} nr.\ 2 (1949), 125--128.
\endthanks
\thanks Received 14 IX 1949
\endthanks
\dedicatory
(Communicated by academician A. N. Kolmogorov, September 14, 1949)
\enddedicatory
\translator
W. Van Assche
\endtranslator
\endtopmatter

\document
In the paper \cite{1}, which deals with the theory of the representation of
Hermitian operators with deficiency index $(p,p)$, we worked out as an
illustration the application to the matrix moment problem.

The present paper includes several important supplements for the investigation
of this problem. These supplements arise in an attempt to construct a complete
and self-contained theory of infinite regular $J_p$-matrices $(p > 1)$ and, just
as  $J_1$-matrices are related to the classical moment problem, this theory of
$J_p$-matrices is related to the matrix moment problem. For clarity, this
analogy led to a natural matrix notation, which considerably facilitates
the subsequent analogy and facts from the theory of $J_p$-matrices and the
matrix moment problem.

We say that the infinite Hermitian matrix $\boldkey A =
\left(a_{i,k}\right)_0^\infty$
is a {\it regular $J_p$-matrix} when it can be written as
$\boldkey A=\left(A_{i,k}\right)_0^\infty$, where $A_{i,k}$ $(i,k=0,1,2,\ldots)$
are square $p\times p$ matrices for which $A_{i,k} = 0$ if $|i-k|>1$ and
$A_{i,i+1}$ $(i=0,1,2,\ldots)$ are nonsingular. In particular,
the Hermitian matrix $\boldkey A=\left(a_{i,k}\right)_0^\infty$ is a regular
$J_p$-matrix
when $a_{i,k}=0$ for $|i-k|>1$\footnote{translator's note: presumably $|i-k|>p$
is meant} and $a_{i,i+p} \neq 0$ for $i=0,1,2,\ldots$.

Our theory of regular $J_p$-matrices is in our opinion of interest since it may
serve as an algebraic model for one-dimensional boundary value problems with one
singularity at an endpoint in a space of functions or vector functions.

Finally we note that infinite regular $J_p$-matrices for $p
> 1$ have already been studied by H. Nagel \cite{2}. However the deepest facts
of the theory of this class of matrices have not been investigated by him.

\subhead 1 \endsubhead
Let $\boldkey A = \left(A_{i,k}\right)_0^\infty$ be a regular $J_p$-matrix
$(p>1)$. Denote by ${\eufm M}_p$ the set of all square $p\times p$ matrices with
complex entries and by ${\eufm B}$ the collection of all polynomials
$P(\lambda)$ given by
$$   P(\lambda) = C_0+ C_1\lambda + \cdots + C_n\lambda^n, $$
where $C_i \in {\eufm M}_p$ $(i=0,1,\ldots,n)$ and $n$ may take any value
from $0,1,2,\ldots$.

We associate with the matrix $\boldkey A$ the sequence of polynomials
$D_k(\lambda)
\in {\eufm B}$ $(k=0,1,2,\ldots)$ defined by the recurrence relation
$$  A_{k,k-1}D_{k-1}(\lambda) + (A_{k,k}-\lambda I) D_k(\lambda)
    + A_{k,k+1} D_{k+1}(\lambda) = 0 \tag 1 $$
$$   (k=0,1,\ldots; D_{-1} \equiv 0), $$
with the polynomial $D_0(\lambda)$ an arbitrary constant nonsingular matrix.

The condition that all the matrices $A_{i,i+1}$ are nonsingular
implies that, for $D_0$ given, all matrix polynomials $D_k(\lambda)$
$(k=1,2,\ldots)$ are uniquely defined and moreover the polynomial $D_k(\lambda)$
has exact degree $k$ with a nondegenerate matrix coefficient for $\lambda^k$
$(k=0,1,2,\ldots)$.

It is easy to see that for each complex number $z$ the limit matrix
$$ H(z) = \lim_{n \to \infty}  \left( \sum_{k=0}^n D_k^*(\bar z)
   D_k(z)\right)^{-1} $$
exists.

Here and in what follows $P^*(\lambda)$ denotes the polynomial obtained from
$P(\lambda) \in {\eufm B}$ by replacing each of its matrix coefficients by its
Hermitian conjugate, so that
$$   P^*(\bar z) = [P(z)]^* . $$

The following results holds (e.g., \cite{2}).
\proclaim{Theorem 1}
The rank $r(z)$ of the Hermitian matrix $H(z)$ is the same for each $z$
belonging to the same half-plane $\text{\rom{Im} } z > 0$ or $\text{\rom{Im} } z
< 0$. \endproclaim

We denote by $\nu_+$ and $\nu_-$ the value of the rank $r(z)$ corresponding
to respectively $\text{\rom{Im} } z > 0$ and $\text{\rom{Im} } z < 0$. The
$J_p$-matrix $\boldkey A$, just as a $Q$-matrix, corresponds to an Hermitian
operator (which we denote by the same letter $\boldkey A$) in a Hilbert space
$\ell^2$ consisting of $x = \{ \xi_k \}_0^\infty$ with complex numbers
$\xi_k$ which are absolutely square summable.

Theorem 1 immediately leads to the following:
\proclaim{Theorem 2}
The numbers $\nu_+$ and $\nu_-$ correspond to the upper and the lower
deficiency index of the Hermitian operator $\boldkey A$.
\endproclaim

In the case of a real matrix $\boldkey A$ it is clear that $\nu_+ = \nu_-$.

\subhead 2 \endsubhead
Every polynomial $P(\lambda) \in {\eufm B}$ of degree  $n$ can be written as
$$   P(\lambda) = \sum_0^n U_kD_k(\lambda), \tag 2 $$
where $U_k \in {\eufm M}_p$ $(k=1,2,\ldots,n)$\footnote{translator's note: $k=0$
should of course also be included}.

For any $P,Q \in {\eufm B}$, where $P$ is given by \thetag{2} and
$$    Q(\lambda) = \sum_0^m V_kD_k(\lambda), $$
we put
$$   \lh P,Q \rh = \sum_0^s U_kV_k^*, \qquad (s=\min(n,m)). $$
In particular
$$   \lh D_i,D_j \rh = \delta_{i,k} I \qquad (i,k=0,1,2,\ldots).
                                                            \tag 3 $$
The ``form'' $\lh P,Q \rh$ is completely defined by condition
\thetag{3} and the following properties:
$$ \lh P_1+P_2,Q \rh = \lh P_1,Q \rh + \lh P_2,Q \rh; \quad
   \lh P,Q_1+Q_2 \rh = \lh P,Q_1 \rh + \lh P,Q_2 \rh; $$
$$ \lh CP,Q \rh = C \lh P,Q \rh; \quad \lh P,CQ \rh = \lh P,Q \rh C^*, \tag 4 $$
where $C$ is any matrix in ${\eufm M}_p$.

It is easy to verify and well-known that by using the relations \thetag{3}
the form $\lh P,Q \rh$ also has the following property:
$$  \lh \lambda P,Q \rh = \lh P, \lambda Q \rh \qquad (P,Q \in {\eufm B}). $$

Let us now form the sequence of matrices
$$   S_n = \lh \lambda^n I,I \rh \qquad (n=0,1,2,\ldots). $$

For each $X_j \in {\eufm M}_p$ $(j=0,1,2,\ldots)$ one has
$$   \sum_{j,k=0}^n  X_j S_{j+k}X_k^* =
     \lh \sum_0^n X_j\lambda^j , \sum_0^n X_j \lambda^j \rh . $$

On the other hand, for each $P \in {\eufm B}$ the expression $\lh P,P \rh$
is an Hermitian matrix (different from zero if $P(\lambda) \not\equiv 0$)
corresponding to a non-negative definite form.
Therefore by choosing matrices $X_j$ $(j=0,1,2,\ldots)$ for which all the
rows, except the first one, consist of only zeros, we see that for arbitrary
$p$-dimensional vectors $x_j = (\xi_{j,1},\xi_{j_2},\ldots,\xi_{j,p})$, not
identically zero, one has
$$   \sum_{j,k=0}^n x_jS_{j+k}x_k^* > 0 \qquad (n=0,1,2,\ldots). \tag 5 $$

It is easy to show that conversely if some sequence of matrices
$\{S_n\}_0^\infty \subset {\eufm M}_p$ satisfies the condition \thetag{5}
then it is induced by some regular $J_p$-matrix.

On the other hand, as has already been shown earlier \cite{1;3}, condition
\thetag{5} is a necessary and sufficient condition for the solvability of the
matrix moment problem
$$   S_n = \int_{-\infty}^\infty \lambda^n\, dT(\lambda) \qquad
(n=0,1,2,\ldots),   \tag 6 $$
where it is required to find an Hermitian matrix function $T(\lambda)$ subject
to the condition that for each $x = (\xi_1,\xi_2,\ldots,\xi_p) \neq 0$ the form
$xT(\lambda)x^*$ is a nondecreasing function of $\lambda \in (-\infty,\infty)$
with an infinite number of points of increase.

In this way every regular $J_p$-matrix corresponds to some matrix moment
problem. Our previous investigation \cite{1} of the problem \thetag{6} makes it
possible, in particular, to state the following:

\proclaim{Theorem 3}
The moment problem \thetag{6} has a unique normalized solution $T(\lambda)$
if and only if one of the numbers $\nu_+$ or $\nu_-$ is equal to zero.
\endproclaim

We say that a solution $T(\lambda)$ of the moment problem is {\it normalized}
if
$$  \lim_{\lambda \to -\infty} T(\lambda) = 0, \qquad
   T(\lambda-0) = T(\lambda)\quad (-\infty < \lambda < \infty). $$

By the norm $\|C\|$ of a matrix $C \in {\eufm M}_p$ we mean the smallest number
$\mu \geq 0$ with the property $C^*C \leq \mu^2 I$ (the inequality for Hermitian
matrices is to be understood as an inequality for the corresponding forms).
\proclaim{Theorem 4}
If $\nu_+=\nu_-=p$ (the case of a completely indeterminate moment problem
\thetag{6}), then the series
$$  H^{-1}(z) = \sum_0^\infty D_k^*(\bar z)D_k(z) $$
converges uniformly on each bounded set of the complex plane
and\footnote{translator's note: presumably the limit for $\|z\| \to \infty$ is
meant}
$$ \lim_{|\lambda| \to \infty} \frac{\log \|H^{-1}(z)\|}{\|z\|} = 0 . $$
Every solution $T(\lambda)$ of the moment problem \thetag{6} satisfies the
inequality
$$   T(\xi+0) - T(\xi-0) \leq H^{-1}(\xi)\qquad (-\infty < \xi < \infty). $$
\endproclaim

For a fixed $\xi$ the equality sign is attained for one and only one normalized
solution $T(\lambda) = T_\xi(\lambda)$ which is determined by the relation
$$ \multline
(z-\xi)^{-1} \left[ I + (z-\xi) \sum_1^\infty E_k^*(z)D_k(\xi) \right]
 \left[\sum_0^\infty D_k^*(z)D_k(\xi) \right]^{-1} \\
 = \int_{-\infty}^\infty \frac{dT_\xi(\lambda)}{\lambda-z} \qquad
(\text{\rom{Im} } z < 0).
\endmultline $$
Here
$$   E_k(z) = \lh \frac{D_k(\lambda)-D_k(z)}{\lambda-z}, I\rh =
  \int_{-\infty}^\infty \frac{D_k(\lambda)-D_k(z)}{\lambda-z}\, dT(\lambda). $$

If $\nu_+=\nu_-=p$, then the series
$$   \sum_1^\infty E_k^*(z)D_k(\zeta), \quad
    \sum_1^\infty E_k^*(z)E_k(\zeta) $$
converge uniformly in the variables $z$ and $\zeta$ on every bounded set of the
complex plane.

Let us construct the entire matrix functions
$$  F_1(z) = I + z \sum_1^\infty  E_k^*(z)D_k(0), \qquad
    F_2(z) =  z \sum_1^\infty  E_k^*(z)E_k(0), $$
$$  G_1(z) = -z \sum_0^\infty  D_k^*(z)D_k(0), \qquad
    G_2(z) = I - z \sum_1^\infty  D_k^*(z)E_k(0). $$
\proclaim{Theorem 5}
If $\nu_+=\nu_-=p$, then each normalized solution $T(\lambda)$ of the problem
\thetag{6} is in a one-to-one correspondence with the collection of all
holomorphic $p\times p$-matrix functions in the
upper half-plane $V(z)$ with $\|V(z)\| \leq 1$ (\/$\text{\rom{Im} } z > 0$),
such that
$$\multline
 \left[ F_1(z)\left(I+V(z)\right) + i F_2(z)\left(I-V(z)\right)\right]
    \left[G_1(z)\left( I+V(z) \right) \right.\\
  \left. + iG_2(z) \left(I-V(z)\right)\right]^{-1} = \int_{-\infty}^\infty
    \frac{dT(\lambda)}{\lambda-z} .
  \endmultline  $$
\endproclaim

We do not turn our attention here to the connection which exists between the
(usual and generalized) resolvent of the operator $\boldkey A$ and the solution
$T(\lambda)$ of problem \thetag{6}, which essentially has been explained in
our paper \cite{1}. We only note that because of this connection each
self-adjoint extension $\tilde{\boldkey A}$ of the operator $\boldkey A$
corresponds to a certain unitary matrix $U$ (and vice versa) in such a way that
the spectrum of the operator $\tilde{\boldkey A}$ coincides with the set of
roots of the equation $\det \left[ G_1(z)(I+U) + iG_2(z)(I-U)\right] = 0$.


\Refs
\ref \no 1 \by M. Krein
\paper Fundamental aspects of the representation theory of Hermitian operators
with deficiency index $(m,m)$
\jour Ukrain. Mat. \v Z. \vol 1 \yr 1949 \pages 3--66
\transl\nofrills English transl. in \jour Amer. Math. Soc. Transl.
(2) \vol 97 \yr 1970 \pages 75--143
\endref
\ref \no 2
\by H. Nagel
\paper \"Uber die aus quadrierbaren Hermiteschen Matrizen entstehenden
Operatoren
\jour Math. Ann. \vol 112 109 \yr 1936 \pages 247--285
\endref
\ref \no 3
\by M. Krein and M. Krasnoselski\u \i
\paper Fundamental theorems on the extension of Hermitian operators and some
 applications to the theory of orthogonal polynomials and the moment problem
\jour Uspehi Mat. Nauk \vol 2  \yr 1947 \pages nr.\ 3 (19), 60--106
\endref
\endRefs
\enddocument